\documentclass[a4paper,14pt]{article}
\usepackage[14pt]{extsizes}
\usepackage{mathtext}                    
\usepackage{cmap}                        
\usepackage[cp1251]{inputenc}            
\usepackage[english, russian]{babel}     
\usepackage[left=2.5cm,right=1cm,top=1cm,bottom=2cm]{geometry} 
\usepackage{amssymb}
\usepackage{amsmath}

\newtheorem{definition}{Definition}
\newtheorem{corollary}{Corollary}
\newtheorem{remark}{Remark}
\newtheorem{lemma}{Lemma}
\newtheorem{example}{Example}

\newtheorem{theorem}{Теорема}

\begin{document}

\centerline{\textbf{Bounded perturbations of Bernstein functions of several operator variables}}


\centerline{A. R. Mirotin}

\centerline{ Department of mathematics, F. Skorina Gomel State University,}

\centerline{ Gomel, 246019, Belarus}

\centerline{amirotin@yandex.ru}


        The paper deals with (multidimensional and one-dimensional) Bochner-Phillips functional calculus.  Bounded perturbations of Bernstein functions of (one or several commuting) semigroup generators on Banach spaces are considered,  conditions for  Lipschitzness and Frechet-differentiability of such functions are obtained,  estimates for the norm of commutators are proved, and a  generalization of Livschits-Kre\u{\i}n  trace formula derived.

\section{Introduction}
\label{1}

The study of the problem of differentiability of functions of self-adjoint operators on
Hilbert space was initiated By Yu. Daletski\u{\i} and S.G. Kre\u{\i}n  \cite{DK}.  Much work has been done during  last decades on the theory of  classes of operator Lipschitz and Frechet-differentiable  functions of (self-adjoint, unitary, or normal) operators on Hilbert space by Birman and Solomyak, Davies, Farforovskaya, Johnson and Williams, Peller, Aleksandrov, Nazarov, Arazy, Barton, Friedman, Pedersen, Shulman, Sukochev, Kissin, Potapov,  Naboko and others.
We refer to \cite{Pel09}, \cite{Pel16}, \cite{ANPell}, \cite{Shul}, \cite{KS2005}, \cite{KS}, \cite{KPSS},  \cite{APPS} as well as to references quoted there for motivation and bibliography. It should be stressed that all these work deal with  Hilbert spaces only. The case of Banach spaces was considered in \cite{RST}.

This paper   is devoted to problems in perturbation theory that arise in an attempt to understand the behavior of the Bernstein function $\psi(A)$ of a semigroup generator $A$ under perturbations of $A$. We consider Bernstein functions of several commuting semigroup generators on Banach spaces (they constitute the subject matter of the so called multidimensional Bochner-Phillips functional calculus).
 We give inter alia conditions for their Lipschitzness and show that  such functions are  $\mathcal{J}$ perturbations preserving where $\mathcal{J}$ is an arbitrary operator ideal; estimates for the norm of commutators are also obtained. In the one-dimensional case   Frechet-differentiability and a  trace formula are  proved. So this work could be considered as a contribution to Bochner-Phillips functional calculus.

One-dimensional Bochner-Phillips functional calculus is a substantial part of the theory of
operator semigroups (see, e.g., \cite{Kish}, \cite{Berg}, \cite{Car}, \cite{MirSMZ}, \cite{SSV}, \cite{MirSF}) and finds important applications in the theory of random
processes (see \cite[ Chap. XIII, Sec. 7]{Fel}, \cite[ Chap. XXIII, Sec. 5]{HiF}, \cite[ Chap. 13, 14]{SSV}, \cite[Chapter 6]{Sato}, and \cite{App}). The foundations of multidimensional
calculus were laid by the author in \cite{Mir97}, \cite{Mir98}, \cite{Mir99}, \cite{SMZ2011},  \cite{IZV2015} (see also \cite{BKM} where the case of multiparametric $C_0$-groups is considered).  Below we recall some notions and
facts from \cite{Boch}, \cite{Mir97}, and \cite{Mir99}, which we need for formulating our results.

\begin{definition}
 \cite{Boch} We say that a nonpositive function $\psi\in C^\infty((-\infty;0)^n)$ {\it
 belongs to the class
 ${\cal T}_n$ } (or is a nonpositive Bernstein function of $n$ variables) if all of its first partial derivatives are absolutely
monotone (a function in $C^\infty((-\infty;0)^n)$ is said to be absolutely monotone if it is nonnegative
together with its partial derivatives of all orders).
\end{definition}

Obviously,  $\psi\in{\cal T}_n$ if and only if $-\psi(-s)$ is a nonnegative Bernstein function of $n$ variables on $(0,\infty)^n$, and ${\cal T}_n$  is a cone under the pointwise addition of functions and multiplication by scalars.
As is known \cite{Boch} (see also \cite{Mir99}, \cite{mz}), each function  $\psi\in {\cal T}_n$  admits an integral representation of the form (here
and in what follows, the dot  denotes inner product in  $\Bbb{R}^n$ and the expression $s\to
-0$ means that $s_1\to -0, \ldots, s_n\to -0$)
$$
\psi(s)=c_0+c_1\cdot s+\int\limits_{\Bbb{R}_+^n\setminus \{0\}}
(e^{s\cdot u}-1)d\mu (u)\quad   (s\in(-\infty;0)^n),  \eqno(\ast)
$$
where $c_0=\psi(-0):=\lim_{s\to -0}\psi(s)$, $c_1=(c_1^j)_{j=1}^n\in \Bbb{R}_+^n$, $c_1^j=\lim_{s_j\downarrow -\infty}\psi(s)/s_j$, and $\mu$ is a positive measure on
$\Bbb{R}_+^n\setminus \{0\}$;  $\mu$
are determined by $\psi$.

A lot of  examples of  Bernstein function of one variable one can fined in \cite{SSV}
(see also \cite{MirSF}, \cite{mz}).

Throughout the paper, $T_{A_1},\dots , T_{A_n}$ denote pairwise commuting one-parameter $C_0$ semigroups
(i.e., strongly continuous semigroups on $\Bbb{R}_+$) on a complex Banach space $X$ with generators $A_1, \dots ,A_n$ respectively satisfying the condition $\|T_{A_j}(t)\|\leq
M_A\quad(t\geq 0, M_A={\rm const}$)  (sometimes we write $T_{j}$ instead of $T_{A_j}$). We denote the domain of $A_j$  by $D(A_j)$ and set $A = (A_1, \dots ,A_n)$. Hereafter, by the commutation of operators $A_1, \dots ,A_n$
we mean the commutation of the corresponding semigroups. By ${\rm Gen}(X)$ we denote the set of all
generators of uniformly bounded $C_0$ semigroups on $X$ and by ${\rm Gen}(X)^n$, the set of all $n$-tuples  $(A_1, \dots ,A_n)$ where $A_j\in {\rm Gen}(X)$. We put $M:=\max\{M_A, M_B\}$ for the pare $A, B$ of $n$-tuples  from  ${\rm Gen}(X)^n$. In the following $\mathcal{L}(X)$ denotes the algebra of linear bounded operators on $X$ and  $I$, the identity operator on $X$.
An operator-valued function $T_A(u) := T_{A_1}(u_1) \dots T_{A_n}(u_n) \quad (u\in
\Bbb{R}_+^n)$ is an $n$-parameter $C_0$ semigroup;
therefore, the linear manifold $D(A):=\cap_{j=1}^nD(A_j)$ is dense in $X$ \cite[ Sec. 10.10]{HiF}.

\begin{definition}
   \cite{Mir99} The value of a function  $\psi\in {\cal T}_n$ of the form ($\ast$) at $A=(A_1,\ldots ,A_n)$ applied
to $x\in D(A)$ is defined by
$$
\psi(A)x=c_0x+c_1\cdot Ax+\int\limits_{\Bbb{R}_+^n\setminus \{0\}}
(T_A(u)-I)xd\mu(u),     
$$
where $c_1\cdot Ax:=\sum_{j=1}^nc_1^jA_jx$.
\end{definition}

Given $\psi\in{\cal T}_n$ and $t\geq 0$, the function $g_t(s):=e^{t\psi(s)}$ is absolutely monotone on $(-\infty;0)^n$. It
is also obvious that $g_t(s)\leq 1$. By virtue of the multidimensional version of the Bernstein-Widder
theorem (see, e.g.,  \cite{Boch}, \cite{BCR}), there exists a unique bounded positive measure $\nu_t$ on $\Bbb{R}_+^n$, such that, for $s\in (-\infty;0)^n$, we have
$$
g_t(s)=\int\limits_{\Bbb{R}_+^n} e^{s\cdot u}d\nu_t(u). 
$$

\begin{definition} In the notation introduced above, we set
$$
g_t(A)x=\int\limits_{\Bbb{R}_+^n} T_A(u)xd\nu_t(u)\ (x\in X) 
$$
(the integral is understood in the sense of Bochner).
\end{definition}

Obviously, $\|g_t(A)\|\leq M_A^n$.  The map $g(A):t\mapsto g_t(A)$ is a $C_0$ semigroup. In the one-dimensional case, it is called the
\textit{semigroup subordinate to} $T_A$ (the terminology is borrowed from probability theory; see \cite[ Sec. X.7]{Fel}
and \cite{App}).
In \cite{Mir09} it was noticed that the closure of the operator $\psi(A)$ exists and is the generator of the $C_0$ semigroup
$g(A)$ (cf. \cite{Mir99}.) It suggests the following final version of the definition of the operator $\psi(A)$.

\begin{definition}
\cite{Mir99} By the value of a function $\psi\in {\cal T}_n$ at an $n$-tuple $A = (A_1,\dots ,A_n)$ of
commuting operators in ${\rm Gen}(X)$ we understand the generator of the semigroup $g(A)$, i.e.,  the closure of the operator defined in the Definition 2. This value is
denoted by $\psi(A)$.
\end{definition}

The functional calculus thus arising is called multidimensional Bochner-Phillips
calculus, or $\mathcal{T}_n$-calculus.

In the sequel unless otherwise stated we  assume for the sake of simplicity  that $c_0=c_1=0$ in the integral representation ($\ast$) of the function $\psi\in \mathcal{T}_n$.

The notation and constraints introduced above are used in what follows without additional
explanations.


\section{Bounded perturbations of Bernstein functions}
\label{2}

\begin{theorem}
(Cf. \cite[Theorem 1]{Naboko}.) Let $\psi\in {\cal T}_n$. Then for every commuting families  $A=(A_1,\dots, A_n)$,  and  $B=(B_1,\dots,B_n)$ from $\mathrm{Gen}(X)^n$ such that the operators $A_i-B_i$ are bounded,  $D(A_i)=D(B_i)\ (i=1,\dots,n)$  the operator $\psi(A)-\psi(B)$ is also bounded and
 $$
 \|\psi(A)-\psi(B)\|\leq -\frac{2e}{e-1}nM^n\psi\left(-\frac{M}{2n}\|A-B\|\right),
 $$
where $\|A-B\|:=(\|A_1-B_1\|,\dots,\|A_n-B_n\|)$.
\end{theorem}

Proof. We have $D(A_i)=D(B_i)\ (i=1,\dots,n)$ and for all $x\in D(A)$
$$
(\psi(A)-\psi(B))x=\int\limits_{\Bbb{R}^n_+\setminus\{0\}} (T_A(u)-T_B(u))xd\mu(u),\eqno(1)
 $$
where $T_A(u):=\prod\limits_{i=1}^nT_{A_i}(u_i)$, $T_B(u):=\prod\limits_{i=1}^nT_{B_i}(u_i)$.

It is easy to verify that
$$
T_A(u)-T_B(u)=\left(\prod\limits_{i=1}^{n-1}T_{A_i}(u_i)\right)(T_{A_n}(u_n)-T_{B_n}(u_n))
$$
$$
+\left(\prod\limits_{i=1}^{n-2}T_{A_i}(u_i)\right)(T_{A_{n-1}}(u_{n-1})-T_{B_{n-1}}(u_{n-1}))T_{B_n}(u_n)
$$
$$
+\dots + T_{A_1}(u_1)(T_{A_2}(u_2)-T_{B_2}(u_2))\prod\limits_{i=3}^nT_{B_i}(u_i)+
(T_{A_1}(u_1)-T_{B_1}(u_1))\prod\limits_{i=2}^nT_{B_i}(u_i).\eqno(2)
$$
It follows that
$$
\|T_A(u)-T_B(u)\|\leq M^{n-1}\sum\limits_{i=1}^{n}\|T_{A_i}(u_i)-T_{B_i}(u_i)\|.\eqno(3)
$$

But  for all $x\in D(A)=D(B)$
$$
(T_{A_i}(u_i)-T_{B_i}(u_i))x=\int\limits_{0}^{u_i}\frac{d}{ds}(T_{B_i}(u_i-s)T_{A_i}(s)x)ds
$$
$$
=\int\limits_{0}^{u_i}T_{B_i}(u_i-s)(A_i-B_i)T_{A_i}(s)xds. \eqno(4)
$$
 Since
 $$
 \|T_{B_i}(u_i-s)(A_i-B_i)T_{A_i}(s)\|\leq M^2\|A_i-B_i\|,
 $$
 both sides in (4) are bounded,  formula (4) holds for all $x\in X$, and
$$
\|T_{A_i}(u_i)-T_{B_i}(u_i)\|\leq M^2\|A_i-B_i\|u_i. \eqno(5)
$$
Now in view of the obvious estimate $\|T_{A_i}(u_i)-T_{B_i}(u_i)\|\leq 2M$, we have
$$
\|T_{A_i}(u_i)-T_{B_i}(u_i)\|\leq 2M\min\{1,(M/2)\|A_i-B_i\|u_i\}.
$$
For reasons of convexity it is clear  that
$$
\min\{1,at\}\leq  \frac{e}{e-1}(1-e^{-at})\ (a, t\geq 0)
$$
 \cite[p. 209]{SSV}. This implies the inequality
$$
\|T_{A_i}(u_i)-T_{B_i}(u_i)\|\leq 2M\frac{e}{e-1}\left(1-e^{-(M/2)\|A_i-B_i\|u_i}\right).
$$
Now we deduce from (3) that
$$
\|T_A(u)-T_B(u)\|\leq \frac{2e}{e-1}M^{n}\sum\limits_{i=1}^{n}\left(1-e^{-(M/2)\|A_i-B_i\|u_i}\right).
$$
Since by the Cauchy inequality
$$
\sum\limits_{i=1}^{n}\left(1-e^{-a_i}\right)\leq n\left(1-e^{-(1/n)\sum_{i=1}^{n}a_i}\right)
$$
$(a_i\geq 0)$,  we have
$$
\|T_A(u)-T_B(u)\|\leq \frac{2e}{e-1}M^{n}n\left(1-e^{-(M/2n)\sum_{i=1}^{n}\|A_i-B_i\|u_i}\right).
$$

Combining this estimate  with the formula (1) we get for all $x\in D(A)$
$$
 \|(\psi(A)-\psi(B))x\|\leq \frac{2e}{e-1}M^{n}n\int\limits_{\Bbb{R}^n_+\setminus\{0\}} \left(1-e^{-(M/2n)\sum_{i=1}^{n}\|A_i-B_i\|u_i}\right)d\mu(u)\|x\|
 $$
 $$
= -\frac{2e}{e-1}M^nn\psi\left(-\frac{M}{2n}\|A-B\|\right)\|x\|,
$$
and the result follows.

\begin{example}
   For every $A, B\in \mathrm{Gen}(X)$ such that $A-B$ is bounded, $D(A)=D(B)$ the following inequalities hold

1) (cf. \cite[formula (10)]{Naboko})
$$
\|(-A)^\alpha-(-B)^\alpha\|\leq \frac{2^{1-\alpha}e}{e-1}M^{1+\alpha}\|A-B\|^\alpha\ (0<\alpha<1);
$$
2)
$$
\|\log(I-A)-\log(I-B)\|\leq \frac{2eM}{e-1}\log\left(1+\frac{M}{2}\|A-B\|\right).
$$
Indeed, the functions $\psi(s)=-(-s)^\alpha (0<\alpha<1)$ and $-\log(1-s)$ belong to $\mathcal{T}_1$.
\end{example}

\begin{corollary}
The multidimensional Bochner-Phillips functional calculus is stable in a sense that
$ \|\psi(A^{(k)})-\psi(B)\|\to 0$  for every $\psi\in {\cal T}_n$ and for every sequence $A^{(k)}$ of commuting families     from $\mathrm{Gen}(X)^n$ such that $\|A^{(k)}-B\|\to 0\ (k\to\infty)$ for some commuting family  $B$   from $\mathrm{Gen}(X)^n$ and $M_{A^{(k)}}\leq M, M= \mathrm{const}$.
\end{corollary}

The next corollary gives (necessary and) sufficient conditions for Bernstein function to be operator Lipschitz  in the class of generators of semigroups of contractions   (the case $M=1$), i.e., in the class of maximally dissipative operators.

\begin{corollary}
Let $\psi\in {\cal T}_n$ be such that $\frac{\partial\psi}{\partial s_i}\left|_{s=-0}\right.\ne \infty$ for all  $i=1,\dots,n$. Then
$$
 \|\psi(A)-\psi(B)\|\leq M^{n+1}\sum_{i=1}^{n}\frac{\partial\psi}{\partial s_i}\left|_{s=-0}\right.\|A_i-B_i\|.
 $$
\end{corollary}

Proof. Using formulae (1), (3), and (5) we have for $x\in D(A)$
 $$
 \|(\psi(A)-\psi(B))x\|\leq \int\limits_{\Bbb{R}^n_+\setminus\{0\}} \|T_A(u)-T_B(u)\|d\mu(u)\|x\|
 $$
 $$
 \leq M^{n-1}\sum_{i=1}^{n}\int\limits_{\Bbb{R}^n_+\setminus\{0\}}
 \|T_{A_i}(u_i)-T_{B_i}(u_i)\|d\mu(u)\|x\|
 $$
 $$
\leq M^{n+1}\sum_{i=1}^{n}\int\limits_{\Bbb{R}^n_+\setminus\{0\}}u_id\mu(u)\|A_i-B_i\|\|x\|
= M^{n+1}\sum_{i=1}^{n}\frac{\partial\psi}{\partial s_i}\left|_{s=-0}\right.\|A_i-B_i\|\|x\|.
 $$

  Recall that a $C_0$ semigroup $T$  is called  exponentially stable if $(\forall t\geq 0)\|T(t)\|\leq Me^{\omega t}$ where $\omega<0$. In the case of generators of exponentially stable semigroups  Corollary 2 can be improved as follows.

\begin{corollary}
Let $\psi\in {\cal T}_n$. For every commuting families  $A=(A_1,\dots, A_n)$  and  $B=(B_1,\dots,B_n)$  from $\mathrm{Gen}(X)^n$ such that the operators $A_i-B_i$ are bounded, $D(A_i)=D(B_i)$, and  semigroups $T_{A_i}$ and $T_{B_i}$ are exponentially stable,  $\|T_{A_i}(t)\|\leq Me^{\omega_i t}$, and $\|T_{B_i}(t)\|\leq Me^{\omega_i t},
\omega_i<0 \quad (i=1,\dots,n)$,  the operator $\psi(A)-\psi(B)$ is also bounded and
 $$
 \|\psi(A)-\psi(B)\|\leq M^{n+1}\sum_{i=1}^{n}\frac{\partial\psi(\omega_i\mathbf{e}_i)}{\partial s_i}\|A_i-B_i\|,
 $$
where $(\mathbf{e}_i)_{i=1}^nХ$ stands for the standard orthogonal basis in $\mathbb{R}^n$.
\end{corollary}

Corollary 3 follows from Corollary 2 applied to the function $\psi(s+\omega):=\psi(s_1+\omega_1,\dots,s_n+\omega_n)$ from $\mathcal{T}_n$ and to the family $A-\omega I:=(A_1-\omega_1 I,\dots,A_n-\omega_n I)$ from $\mathrm{Gen}(X)^n$.

In the case $n=1$ the following theorem (for contraction semigroups) was proved in \cite[Corollary 13.9]{SSV}.

\begin{theorem}
 Let $\psi\in {\cal T}_n$. Then for every commuting families  $A=(A_1,\dots, A_n)$, and  $B=(B_1,\dots,B_n)$ such that   operators $A_i$ and $B_i$ belongs to $\mathrm{Gen}(X)$ and commute  $(i=1,\dots,n)$ the  following estimate holds for all $x\in D(A)\cap D(B)$
 $$
 \|(\psi(A)-\psi(B))x\|\leq -\frac{2e}{e-1}nM^n\psi\left(-\frac{M}{2n}\|(A-B)x\|\right),
 $$
where $\|(A-B)x\|:=(\|(A_1-B_1)x\|,\dots,\|(A_n-B_n)x\|)$.
\end{theorem}

Proof.  Since operators $A_i$ and $B_i$ commute, it follows from (4) that  for all $x\in D(A)\cap D(B)$
 $$
\|(T_{A_i}(u_i)-T_{B_i}(u_i))x\|\leq M^2\|(A_i-B_i)x\|u_i\ (i=1,\dots,n).
$$
 Using this estimate instead of (5) we  can proceed exactly as in the proof of Theorem 1.

\begin{corollary}
(Cf. \cite[Theorem 16]{SMZ2011}.) Let  $\psi\in {\cal T}_n$. For every commuting family  $A=(A_1,\dots, A_n)$  from $\mathrm{Gen}(X)$   the following inequality holds
$$
 \|\psi(A)x\|\leq -\frac{2e}{e-1}nM^n\psi\left(-\frac{M}{2n}\|Ax\|\right),
 $$
where $x\in D(A), \|Ax\|:=(\|A_1x\|,\dots,\|A_nx\|)$.

In particular $\psi(A)$ is bounded if  $\psi$ is bounded.
\end{corollary}

Proof. This corollary follows from Theorem 2 with $B=O$.

\begin{corollary}
 Let $\psi\in {\cal T}_n$ be such that $\frac{\partial\psi}{\partial s_i}\left|_{s=-0}\right.\ne \infty$ for all  $i=1,\dots,n$. If operators $A_i$ and $B_i$ from $\mathrm{Gen}(X)$ commute  $(i=1,\dots,n)$ the  following estimate holds for all $x\in D(A)\cap D(B)$
$$
 \|(\psi(A)-\psi(B))x\|\leq M^{n+1}\sum_{i=1}^{n}\frac{\partial\psi}{\partial s_i}\left|_{s=-0}\right.\|(A_i-B_i)x\|.
 $$
\end{corollary}

This Corollary  is proved similarly to Corollary 2.

\begin{definition}
 \cite{Voigt} A closed (with respect to the operator norm) ideal $\mathcal{E}$ of the algebra $\mathcal{L}(X)$ of bounded operators on $X$ is defined to have the \textit{strong convex compactness property } if for every finite measure space $(\Omega,\mu)$ and every strongly measurable bounded function $U:\Omega\to \mathcal{L}(X)$ the strong integral $\int_\Omega Ud\mu$  belongs to $\mathcal{E}$.
\end{definition}

Theorems 3 and 4 below show in particular that functions from $\mathcal{T}_n$ are $\mathcal{J}$-stable in the sense of \cite{KS2005}.

\begin{theorem}
 Let $\mathcal{E}$ be a closed ideal of the algebra $\mathcal{L}(X)$ which has the strong convex compactness property and  $\psi\in {\cal T}_n$. For every commuting families  $A=(A_1,\dots, A_n)$  and  $B=(B_1,\dots,B_n)$ from $\mathrm{Gen}(X)^n$ such that  $A_i-B_i\in\mathcal{ E}, D(A_i)=D(B_i)\quad (i=1,\dots,n)$  the operator $\psi(A)-\psi(B)$ belongs to $\mathcal{E}$, too.
 \end{theorem}

Proof.  Since  $\mathcal{E}$ is a closed ideal in $\mathcal{L}(X)$, formula (4) implies that $T_{A_i}(u_i)-T_{B_i}(u_i)\in \mathcal{E}$ for all $i$ (the Bochner integral in the right hand side of (4) exists in  the operator norm).
The operator $\psi(A)-\psi(B)$ is bounded by Theorem 1 and therefore  for all $x\in X$
 $$
(\psi(A)-\psi(B))x=\int\limits_{\Bbb{R}^n_+\setminus [0,\delta)^n} +\int\limits_{[0,\delta)^n\setminus\{0\}} (T_A(u)-T_B(u))xd\mu(u),\eqno(6)
 $$
 where $\delta>0$ and both integrals converge in the strong operator topology. Because of the set $\Bbb{R}^n_+\setminus [0,\delta)^n$ has finite $\mu$-measure for  every $\delta>0$ \cite[Lemma 3.1]{Mir99}, and the function $T_A(u)-T_B(u)$ is uniformly bounded, the first summand in the right-hand side of (6) belongs to $\mathcal{E}$ by the strong convex compactness property.

 As regards to the second summand in (6), the formula (2) enables us to present it as the finite sum of strong integrals  of the form
 $$
 \int\limits_{[0,\delta)^n\setminus\{0\}} V_i(u)\frac{1}{\sum_{j=1}^nu_j} (T_{A_i}(u_i)-T_{B_i}(u_i))W_i(u)x\left(\sum_{j=1}^nu_j\right)d\mu(u),\eqno(7)
 $$
 where $V_i(u)$ and $W_i(u)$ are uniformly bounded strongly continuous operator-valued functions.
 Moreover, formula (4) implies also that
 $$
 \frac{1}{\sum_{j=1}^nu_j}\|T_{A_i}(u_i)-T_{B_i}(u_i)\|\leq
M^2\|A_i-B_i\|\frac{u_i}{\sum_{j=1}^nu_j}\leq
M^2\|A_i-B_i\|,
 $$
and so integrands in (7) are bounded and strongly continuous functions in $u\in {\Bbb{R}^n_+\setminus\{0\}}$.
Since by  \cite[Lemma  3.1]{Mir99} the measure  $(\sum_{j=1}^nu_j)d\mu(u)$ is bounded on $[0,\delta)^n$,
all the operators of the form (7) belongs to $\mathcal{E}$ by the strong convex compactness property.

\begin{corollary}
Let  $\psi\in {\cal T}_n$. For every commuting families  $A=(A_1,\dots, A_n)$  and  $B=(B_1,\dots,B_n)$ from $\mathrm{Gen}(X)^n$ such that $D(A_i)=D(B_i)$ and  $A_i-B_i$ are compact $(i=1,\dots,n)$ the operator $\psi(A)-\psi(B)$ is compact, too.
\end{corollary}

Indeed, by \cite[Theprem 1.3]{Voigt} the ideal of compact operators on $X$ possesses the strong convex compactness property.

\begin{remark} See \cite{Voigt} for other examples of operator ideals  with strong convex compactness property.
\end{remark}

\begin{corollary}
 (Cf. \cite[Theorem 1]{Naboko}.)  Let  $\psi\in {\cal T}_n$. For every commuting family  $A=(A_1,\dots, A_n)$ of compact operators from $\mathrm{Gen}(X)$  the operator $\psi(A)$ is compact.
\end{corollary}

 Corollary 7 follows from  Corollary 6 (we assumed at the beginning that
 $\psi(0)=0$).

Below we shall assume that  a (two sided) operator ideal $(\mathcal{J}, \|\cdot\|_{_\mathcal{J}})$ on $X$ is \textit{symmetrically normed} in the sense that $\|ASB\|_{_\mathcal{J}}\leq \|A\|\|S\|_{_\mathcal{J}}\|B\|$
for $A, B\in\mathcal{ L}(X)$ and $S\in \mathcal{J}$.  The following theorem shows that Bernstein functions are  $\mathcal{J}$ perturbations preserving.

\begin{theorem}
Let $(\mathcal{J}, \|\cdot\|_{_\mathcal{J}})$ be an operator ideal on $X$ and  $\psi\in {\cal T}_n$ be such that $\frac{\partial\psi}{\partial s_i}\left|_{s=-0}\right.\ne \infty$ for all  $i=1,\dots,n$. For every commuting families  $A=(A_1,\dots, A_n)$  and  $B=(B_1,\dots,B_n)$ from $\mathrm{Gen}(X)^n$ such that  $A_i-B_i\in\mathcal{ J}, D(A_i)=D(B_i)\ (i=1,\dots,n)$  the operator $\psi(A)-\psi(B)$ belongs to $\mathcal{J}$, too, and
$$
 \|\psi(A)-\psi(B)\|_{_\mathcal{J}}\leq M^{n+1}\sum_{i=1}^{n}\frac{\partial\psi}{\partial s_i}\left|_{s=-0}\right.\|A_i-B_i\|_{_\mathcal{J}}.
 $$
\end{theorem}

Proof. Since $T_{B_i}(u_i-s)(A_i-B_i)T_{A_i}(s)\in \mathcal{J}$  and
$$
\|T_{B_i}(u_i-s)(A_i-B_i)T_{A_i}(s)\|_{_\mathcal{J}}\leq M^2\|A_i-B_i\|_{_\mathcal{J}}
$$
for all $s$, $0\leq s\leq u_i$  $(i=1,\dots,n)$, formulae (4) and (2) entail that $T_{B}(u)-T_{A}(u)\in \mathcal{J}$ and
$$
 \|T_{B}(u)-T_{A}(u)\|_{_\mathcal{J}}\leq M^{n+1}\sum_{i=1}^{n}\|A_i-B_i\|_{_\mathcal{J}}u_i
 $$
for all $u\in\Bbb{ R}_+^n$. It follows that
$$
 \int\limits_{\Bbb{ R}_+^n\setminus\{0\}}\|T_{B}(u)-T_{A}(u)\|_{_\mathcal{J}}d\mu(u)\leq M^{n+1}\sum_{i=1}^{n}\|A_i-B_i\|_{_\mathcal{J}}\int\limits_{\Bbb{ R}_+^n\setminus\{0\}}u_id\mu(u)
 $$
 $$
= M^{n+1}\sum_{i=1}^{n}\frac{\partial\psi}{\partial s_i}\left|_{s=-0}\right.\|A_i-B_i\|_{_\mathcal{J}}.
 $$
In particular  the operator in  the right hand side of the formula (1) belongs to $\mathcal{J}$. Since by Theorem 1 the operator in  the left hand side of  this formula is bounded,   (1) holds for all $x\in X$, and the proof is complete.

Arguing as in the proof of Corollary 3, we get also for generators of exponentially stable semigroups the following

\begin{corollary}
Let $(\mathcal{J}, \|\cdot\|_{_\mathcal{J}})$ be an operator ideal on $X$ and  $\psi\in {\cal T}_n$. For every commuting families  $A=(A_1,\dots, A_n)$  and  $B=(B_1,\dots,B_n)$ from $\mathrm{\mathrm{Gen}}(X)^n$ such that  $A_i-B_i\in\mathcal{ J}, D(A_i)=D(B_i)$, and $\|T_{A_i}(t)\|\leq Me^{\omega_i t}$,  $\|T_{B_i}(t)\|\leq Me^{\omega_i t},
\omega_i<0, \ (i=1,\dots,n)$, the operator $\psi(A)-\psi(B)$ belongs to $\mathcal{J}$, too, and
$$
 \|\psi(A)-\psi(B)\|_{_\mathcal{J}}\leq M^{n+1}\sum_{i=1}^{n}\frac{\partial\psi(\omega_i\mathbf{e}_i)}{\partial s_i}\|A_i-B_i\|_{_\mathcal{J}}.
 $$
\end{corollary}

\section{Estimates for the norm of commutators}
\label{3}

First recall that an operator $H\in \mathcal{L}(X)$ is called \textit{Hermitian } if $\|e^{isH}\|=1$ for all $s\in \mathbb{R}$. We put $V_H(s)=e^{isH}$ for short; this is an automorphism of $X$  for every $s\in \mathbb{R}$. We write $[A,B]$ for the commutator $AB-BA$. If this operator is dense defined and bounded we denote by $[A,B]$ its extension to $X$, too.

\begin{lemma}
Let $H\in \mathcal{L}(X)$ be Hermitian, $T_A$ a $C_0$-semigroup on $X$ with generator $A$, and $V_H(s)$ maps $D(A)$ into itself for all $s\in \mathbb{R}$. Then for all $x\in D(A), s\in \mathbb{R}$
$$
[A,V_H(s)]x=is\int\limits_0^1V_H(sr)[A,H]V_H(s(1-r))xdr.
$$
\end{lemma}

Proof. First we prove that
$$
[B,V_H(s)]=is\int\limits_0^1V_H(sr)[B,H]V_H(s(1-r))dr
$$
for all $B\in \mathcal{L}(X)$. Indeed, the integral in the right-hand side of this formula exists in the sense of Bochner with respect to the operator norm and
$$
B-V_H(s)BV_H(s)^{-1}=\int\limits_0^1\frac{d}{dr}\left(-V_H(sr)BV_H(-sr)\right)dr=is\int\limits_0^1V_H(sr)[B,H]V_H(-sr))dr.
$$
Now taking $B=T_A(t)$ we  obtain that for  $x\in D(A)$
$$
\frac{1}{t}[T_A(t),V_H(s)]x=is\int\limits_0^1V_H(sr)\frac{1}{t}[T_A(t),H]V_H(s(1-r))xdr\quad (t>0),
$$
and the result follows as $t\to 0$.

\begin{corollary}
If, in addition, the operator $[A,H]$ belongs to an operator ideal $\mathcal{J}$ on $X$,  the operator $[A,V_H(s)]$ belongs to $\mathcal{J}$, too, and
$$
\|[A,V_H(s)]\|_{_\mathcal{J}}\leq |s|\|[A,H]\|_{_\mathcal{J}}\quad (s\in \mathbb{R}).
$$
\end{corollary}

In fact, if $[A,H]\in \mathcal{J}$ the integral $\int_0^1V_H(sr)[A,H]V_H(s(1-r))dr$ exists in the sense of Bochner with respect to the $\mathcal{J}$-norm for every $s\in \mathbb{R}$.

\begin{theorem}
 (Cf. \cite[Theorem 3.5]{KS}.)   Let $H\in \mathcal{L}(X)$ be Hermitian, $[A_j,H]$ belongs to an operator ideal $\mathcal{J}$ on $X$, and $V_H(s)$ maps $D(A_j)$ into itself for all $s\in \mathbb{R}$  ($j=1,\dots, n$). Then for every $\psi\in {\cal T}_n$  such that $\frac{\partial\psi}{\partial s_j}\left|_{s=-0}\right.\ne \infty$   ($j=1,\dots,n$) and $\psi(A)$ belongs to $\mathcal{J}$, the following inequality holds
$$
\|[\psi(A),H]\|_{_\mathcal{J}}\leq  M^{n+1}\sum_{j=1}^{n}\frac{\partial\psi}{\partial s_j}\left|_{s=-0}\right.\|[A_j,H]\|_{_\mathcal{J}}.
$$
\end{theorem}

Proof. Let
$$
V_H(s)AV_H(s)^{-1}:=(V_H(s)A_1V_H(s)^{-1},\dots,V_H(s)A_nV_H(s)^{-1}).
 $$
 Note that  $V_H(s)A_jV_H(s)^{-1}\in \mathrm{Gen}(X)$ and
$V_H(s)T_{A}V_H(s)^{-1}=T_{V_H(s)AV_H(s)^{-1}}$  ($s\in \mathbb{R}$). Moreover, for $t>0$ and $x\in X$
we have
$$
g_t(V_H(s)AV_H(s)^{-1})x=\int\limits_{\mathbb{R}_+^n}(T_{V_H(s)AV_H(s)^{-1}}(u)-I)xd\nu_t(u)=V_H(s)g_t(A)V_H(s)^{-1}x.
$$
It follows that $\psi(V_H(s)AV_H(s)^{-1})=V_H(s)\psi(A)V_H(s)^{-1}$. Now applying Theorem 4 and Corollary 9 we obtain
$$
\|[\psi(A),V_H(s)]\|_{_\mathcal{J}}=\|(\psi(A)-\psi(V(s)AV_H(s)^{-1}))V_H(s)\|_{_\mathcal{J}}
$$
$$
\leq \|\psi(A)-\psi(V_H(s)AV_H(s)^{-1})\|_{_\mathcal{J}}
\leq  M^{n+1}\sum_{j=1}^{n}\frac{\partial\psi}{\partial s_j}\left|_{s=-0}\right.\|A_j-V_H(s)A_jV_H(s)^{-1}\|_{_\mathcal{J}}
$$
$$
\leq M^{n+1}\sum_{j=1}^{n}\frac{\partial\psi}{\partial s_j}\left|_{s=-0}\right.\|[A_j,V_H(s)]\|_{_\mathcal{J}}
\leq |s|M^{n+1}\sum_{j=1}^{n}\frac{\partial\psi}{\partial s_j}\left|_{s=-0}\right.\|[A_j,H]\|_{_\mathcal{J}}.
$$
 Dividing by $|s|$, we obtain the desired inequality.

\begin{corollary}
Let operators $[A_j,H]$ are bounded, the function $\psi\in\mathcal{ T}_n$ is bounded or each $A_j$ is bounded,
$\frac{\partial\psi}{\partial s_j}\left|_{s=-0}\right.\ne \infty$, and $V_H(s)$ maps $D(A_j)$ into itself for all $s\in \mathbb{R}$   ($j=1,\dots,n$). Then
$$
\|[\psi(A),H]\|\leq  M^{n+1}\sum_{j=1}^{n}\frac{\partial\psi}{\partial s_j}\left|_{s=-0}\right.\|[A_j,H]\|.
$$
\end{corollary}

Indeed, since $\psi(A)$ is bounded (see Theorem 1 and Corollary  4), one can apply Theorem 5 to the operator ideal $\mathcal{L}(X)$.

In the case of  exponentially stable semigroups  Theorem 5 can be improved, as well.

In the following we consider the case $n=1$.

\section{Differentiability}
\label{4}

\begin{definition}
  (Cf. \cite{KPSS}.) Let $(\mathcal{J}, \|\cdot\|_{_\mathcal{J}})$ be an operator ideal on $X$, $\psi\in {\cal T}_1$, $A\in \mathrm{Gen}(X)$. We call the bounded linear operator   $\psi_A^\nabla$ on $\mathcal{J}$ (transformator)   the \textit{$\mathcal{J}$-Frechet derivative} of the operator function $\psi:\mathrm{Gen}(X)\to \mathrm{Gen}(X)$  at  the point $A$, if for every operator $\Delta A\in \mathcal{J}$ such that $A+\Delta A\in \mathrm{Gen}(X)$ we have
$$
\|\psi(A+\Delta A)-\psi(A)-\psi_A^\nabla(\Delta A)\|_{_\mathcal{J}}=o(\|\Delta A\|_{_\mathcal{J}}) \mbox{ as }  \|\Delta A\|_{_\mathcal{J}}\to 0.
$$
\end{definition}

Evidently, the Frechet derivative at    the point $A$ is unique.

Before we formulate the next theorem recall that  if $\psi^\prime(-0)\ne\infty$ the derivative  $\psi^\prime (s)$ of a function  $\psi\in {\cal T}_1$ equals to  $\int_{(0,\infty)}e^{sv}vd\mu(v) (s\leq 0)$ and absolutely monotonic on $(-\infty,0]$. So for every  $A\in \mathrm{Gen}(X)$ the operator
$$
\psi^\prime(A):=\int\limits_{(0,\infty)}T_A(v)vd\mu(v)
$$
(here $T_A$ denotes the $C_0$-semigroup generated by $A$) exists in the sense of the Hille-Phillips functional calculus \cite[Definition 15.2.2]{HiF}
and belongs to $\mathcal{L}(X)$.

\begin{theorem}
Let $\psi\in {\cal T}_1$, $\psi^\prime(-0)\ne\infty$, $A\in \mathrm{Gen}(X)$.
The $\mathcal{ L}(X)$-Frechet derivative for the operator function $\psi$  at    the point $A$ exists and equals to $\psi^\prime(A)$ in a sense that $\psi_A^\nabla(B)=\psi^\prime(A)B$ for every $B\in\mathcal{ L}(X)$.
\end{theorem}

Proof.
For the proof we need the following

\begin{lemma}
 \cite[Theorem 5]{mz}. Let the function $\psi\in {\cal T}_1$ has the integral representation ($\ast$), and $\psi^\prime(-0)\ne\infty$.
Then the function
$$
 \varphi(s_1, s_2) :=\left\{
\begin{array}
{@{\,}r@{\quad}l@{}}
\frac{\psi(s_1)-\psi(s_2)}{s_1-s_2}-\psi^\prime(-0) \quad {\rm if }\  s_1\ne s_2,\\
\psi^\prime(s_1)-\psi^\prime(-0)  \quad {\rm if } \ s_1= s_2
\end{array}
\right.
$$
belongs to  ${\cal T}_2$ and has the integral representation

$$
\varphi(s_1,s_2)=\int\limits_{\mathbb{R}^2_+\setminus\{0\}}
(e^{s_1u_1+s_2u_2}-1) d\mu_1(u_1,u_2),\eqno(8)
$$
\textit{where $d\mu_1(u_1,u_2)$ is the image of the measure $1/2d\mu(v)dw$ under the mapping}
 $u_1=(v+w)/2, u_2=(v-w)/2$.
\end{lemma}

Now we claim that for every $A_1,A_2\in {\rm Gen}(X)$, such that $A_1-A_2\in \mathcal{L}(X)$ the following equality holds for $x\in D(A_1) = D(A_2)$
$$
\varphi(A_1,A_2)(A_1-A_2)x=(\psi(A_1)-\psi(A_2))x-\psi^\prime(-0)(A_1-A_2)x.\eqno(9)
$$
For the proof first note that
 in view of (8) for $x\in D(A_1)$ we have
$$
\varphi(A_1,A_2)(A_1-A_2)x=\int\limits_{\mathbb{R}^2_+\setminus\{0\}} (T_1(u_1)T_2(u_2)-I)(A_1-A_2)x
d\mu_1(u_1,u_2)
$$
(for simplicity we write $T_i$ instead of  $T_{A_i}, i=1,2)$. Let  $\Omega$ be the angle in the $(v,w)$ plane bounded by the bisectors of the first and fourth quadrants. If we put in the last integral  $v=u_1+u_2, w=u_1-u_2$, then $(v,w)$ runs over $\Omega\setminus\{0\}$ and
 we get
$$
\varphi(A_1,A_2)(A_1-A_2)x= \frac{1}{2}\int\limits_{(0,\infty)} d\mu(v)\int\limits_{-v}^v
\left(T_1\left(\frac{v+w}{2}\right)T_2\left(\frac{v-w}{2}\right)-I\right)(A_1-A_2)xdw
$$
$$
=\int\limits_{(0,\infty)} d\mu(v)\frac{1}{2} \int\limits_{-v}^v
T_1\left(\frac{v+w}{2}\right)T_2\left(\frac{v-w}{2}\right)(A_1-A_2)xdw-\psi^\prime(-0)(A_1-A_2)x.\eqno(10)
$$
Consider the identity ($x\in D(A_1))$
$$
\frac{1}{2} \int\limits_{-v}^v
T_1\left(\frac{v+w}{2}\right)T_2\left(\frac{v-w}{2}\right)(A_1-A_2)xdw
$$
$$
=\frac{1}{2}\int\limits_{-v}^v
T_1\left(\frac{v+w}{2}\right)T_2\left(\frac{v-w}{2}\right)A_1xdw-
\frac{1}{2}\int\limits_{-v}^v
T_2\left(\frac{v-w}{2}\right)T_1\left(\frac{v+w}{2}\right)A_2xdw.\eqno(11)
$$
Because of
$$
\frac{1}{2} \int\limits_{-v}^v
T_1\left(\frac{v+w}{2}\right)T_2\left(\frac{v-w}{2}\right)A_1xdw
$$
$$
=\frac{1}{2} \int\limits_{-v}^v
T_2\left(\frac{v-w}{2}\right)A_1T_1\left(\frac{v+w}{2}\right)xdw
=\int\limits_{-v}^v
T_2\left(\frac{v-w}{2}\right)d_wT_1\left(\frac{v+w}{2}\right)x
$$
$$
=T_2\left(\frac{v-w}{2}\right)T_1\left.\left(\frac{v+w}{2}\right)x\right|_
{w=-v}^{w=v}-
\int\limits_{-v}^v
T_1\left(\frac{v+w}{2}\right)\left(-\frac{1}{2}\right)A_2T_2\left(\frac{v-w}{2}
\right)xdw
$$
$$
=T_1(v)x-T_2(v)x+\frac{1}{2}\int\limits_{-v}^v
T_2\left(\frac{v-w}{2}\right)T_1\left(\frac{v+w}{2}\right) A_2xdw,
$$
the formula (11) implies that
$$
\frac{1}{2} \int\limits_{-v}^v
T_1\left(\frac{v+w}{2}\right)T_2\left(\frac{v-w}{2}\right)(A_1-A_2)xdw=
(T_1(v)-T_2(v))x,
$$
and (9) follows from (10).

Now putting  $A_2=A, A_1-A_2=\Delta A$ in the formula (9)  ($A\in \mathrm{Gen}(X), \Delta A\in \mathcal{L}(X)$), we have for $x\in D(A_1)$
$$
(\psi(A+\Delta A)-\psi(A))x=\varphi(A+\Delta A,A)\Delta Ax+\psi^\prime(-0)\Delta Ax.\eqno(12)
$$
By Theorem 1 the operator
$$
\alpha(\Delta A):=\varphi(A+\Delta A,A)-\varphi(A,A)
$$
is bounded and
$$
\|\alpha(\Delta A)\|\leq -\frac{4eM^2}{e-1}\varphi\left(-\frac{M}{4}\|\Delta A\|,0\right)\to -\frac{4eM^2}{e-1}\varphi(0,0)=0 \mbox{ as } \|\Delta A\|\to 0.
$$
Thus the formula (12)  entails the equality
$$
(\psi(A+\Delta A)-\psi(A))x=\varphi(A,A)\Delta Ax+\psi^\prime(-0)\Delta Ax+\alpha(\Delta A)\Delta Ax, 
$$
and Theorem 6 will be proved if we establish that
$$
\varphi(A,A)=\psi^\prime(A)-\psi^\prime(-0)I.
$$
 To this end note that the Definition 2 implies in view of formula (8) that for $x\in D(A_1)$
$$
\varphi(A,A)x=\int\limits_{\mathbb{R}^2_+\setminus\{0\}}
(T_2(u_1+u_2)-I)x d\mu_1(u_1,u_2).
$$
 If we put here  $v=u_1+u_2, w=u_1-u_2$ as in the proof of the formula (10),  we get
$$
\varphi(A,A)x=\int\limits_{\Omega\setminus\{0\}}(T_2(v)-I)x\frac{1}{2}d\mu(v)dw=\frac{1}{2}\int\limits_{(0,\infty)}d\mu(v)
\int\limits_{-v}^v(T_2(v)-I)xdw
$$
$$
=\int\limits_{(0,\infty)}T_2(v)xvd\mu(v)-\int\limits_{(0,\infty)}vd\mu(v)x=
\psi^\prime(A)x-\psi^\prime(-0)x.
$$
This completes the proof.

Note that the condition $\psi^\prime(-0)\ne\infty$ is also necessary for the Frechet differentiability of the function $\psi$ at every point $A\in \mathrm{Gen}(X)$ (take $A=O$) but the following
corollary holds.

\begin{corollary}
Let $\psi\in {\cal T}_1$,  $A\in \mathrm{Gen}(X)$, and $\|T_{A}(t)\|\leq Me^{\omega t}$,
$\omega<0$.
Then the $\mathcal{ L}(X)$-Frechet derivative for the operator function $\psi$ at the point $A$ exists and equals to $\psi^\prime(A)$.
\end{corollary}

Proof.
 To use Theorem 6 we need the condition $\psi'(-0)\ne \infty$. To bypass it we apply Theorem 6 to the function $\psi(s+\omega)$ from $\mathcal{T}_1$ and  to the operator $A-\omega I$ from $\mathrm{Gen}(X)$.

\begin{theorem}
Let $(\mathcal{J}, \|\cdot\|_{\mathcal{J}})$ be an operator ideal on $X$,  $\psi\in {\cal T}_1$, and $\psi^\prime(-0)\ne\infty,\psi^{\prime\prime}(-0)\ne\infty$. For every $A\in \mathrm{Gen}(X)$
the $\mathcal{J}$-Frechet derivative  for the operator function $\psi$ at the point $A$  exists and equals to $\psi^\prime(A)$ in a sense that $\psi_A^\nabla(B)=\psi^\prime(A)B$ for every $B\in \mathcal{J}$.
\end{theorem}

Proof. Indeed, if in the proof of Theorem 6 we put $\Delta A\in \mathcal{J}$, Theorem 4 implies that the operator
$\alpha(\Delta A)$ belongs to $\mathcal{J}$, too and
$$
\|\alpha(\Delta A)\|_{_\mathcal{J}}\leq M^2\frac{\partial\varphi(-0,-0)}{\partial s_1}\|\Delta A\|_{_\mathcal{J}}=M^2\frac{1}{2}\psi^{\prime\prime}(-0)\|\Delta A\|_{_\mathcal{J}},
$$
since, by  Taylor's formula  (below $\xi$ lies between $s_1$ and $s_2$),
$$
\frac{\partial\varphi(-0,-0)}{\partial s_1}
=\lim_{(s_1, s_2)\to(-0-0)}\frac{\partial\varphi(s_1,s_2)}{\partial s_1}
=\lim_{(s_1, s_2)\to(-0,-0)}\frac{\psi^\prime(s_1)(s_1-s_2)-(\psi(s_1)-\psi(s_2))}{(s_1-s_2)^2}
$$
$$
=\lim_{(s_1, s_2)\to(-0,-0)}\frac{1}{2}\psi^{\prime\prime}(\xi)=\frac{1}{2}\psi^{\prime\prime}(-0).
$$

The remaining part of the proof is the same as in  Theorem 6.

In context of Theorem 7 there is an analog of Corollary 11  for exponentially stable semigroups, as well.

\begin{corollary}
Let $(\mathcal{J}, \|\cdot\|_{_\mathcal{J}})$ be an operator ideal on $X$,  $\psi\in {\cal T}_1$. Let $A\in \mathrm{Gen}(X)$, and $\|T_{A}(t)\|\leq Me^{\omega t}$,
$\omega<0$.
Then the  $\mathcal{J}$-Frechet derivative  for the operator function $\psi$ at the point $A$ exists and equals to $\psi^\prime(A)$.
\end{corollary}

\begin{remark} The map $A\mapsto \psi^\prime(A)$ is analytic in the sense of Hille and Phillips \cite[Theorem 15.5.2]{HiF}.
\end{remark}

\section{Trace formula}
\label{5}

The trace formula for a trace class perturbation of a self-adjoint operator  was
proved in a special case in \cite{L} (where its physical applications were also discussed) and in the general case in
 \cite{Kr1}. A survey of farther developments (in  context of  Hilbert spaces) and bibliography one can fined in \cite{BY} and \cite{Pel09}.

 In this section we introduce a spectral shift function and prove  a Livschits-Kre\u{\i}n trace formula for a trace class perturbation of a generator of bounded holomorphic semigroup in Banach space. Recall  that if the Banach space $X$ has the approximation property there is a continuous linear functional $\mathrm{tr}$ of norm 1  (\textit{a trace}) on the operator ideal $(\mathfrak{S}_1, \|\cdot \|_{\mathfrak{S}_1})$ of nuclear operators on $X$ (see, e. g., \cite[p. 64]{DF}).

\begin{theorem}
Let the Banach space $X$ has the approximation property.  Let  $A$ and $B$ be generators of  $C_0$-semigroups $T_A$ and $T_B$ respectively on $X$ holomorphic in the  half plane $\{\mathrm{Re}(z)>0\}$ and
  satisfying $\|T_A(z)\|, \|T_B(z)\|\leq M$  $(\mathrm{Re}(z)>0)$. If $A-B\in\mathfrak{ S}_1$ there exists a unique tempered distribution $\xi$ supported in $\mathbb{R}_+$  such that for every $\psi\in \mathcal{T}_1$ with $\psi^\prime(-0)\ne\infty$ we have
$$
\mathrm{tr}(\psi(A)-\psi(B))=\int\limits_{(0,\infty)}\langle\xi(t),e^{-ut}\rangle ud\mu(u),
$$
where  (as above)  $\mu$ stands for the  representing measure of $\psi$ and  $\langle\xi(t),\varphi(t)\rangle$ denotes the action of a distribution on a test function.
\end{theorem}

Proof. Consider the function
$$
F(z)=\frac{1}{z}(T_A(z)-T_B(z))\quad (\mathrm{Re}(z)>0).
$$
Theorem 4 implies that  $F:\{\mathrm{Re}(z)>0\}\to \mathfrak{S}_1$. Arguing as in the proof of the formula (4) we have
$$
F(z)=\frac{1}{z}\int\limits_{[0,z]}T_{B}(z-s)(A-B)T_{A}(s)ds\eqno(13)
$$
where the integral exists in the sense of Bochner  with respect to the trace norm (and with respect to the operator  norm) because the norm of the integrand in (13) does not exceed
$M^2\|A-B\|_{\mathfrak{S}_1}$. It follows that $F$ is bounded, since for each $z\in \mathbb{C}$ with $\mathrm{Re}(z)>0$
$$
\|F(z)\|_{\mathfrak{S}_1}\leq \frac{1}{|z|}\max_{s\in[0,z]}\|T_{B}(z-s)(A-B)T_{A}(s)\|_{\mathfrak{S}_1}|z| \leq M^2\|A-B\|_{\mathfrak{S}_1}.
$$
We clame that $F$ is continuous. In fact, let $g(z):=T_A(z)-T_B(z)$. For every complex $z$ with  $\mathrm{Re}(z)>0$ and sufficiently small $\Delta z$ we have in view of (13)
$$
g(z+\Delta z)-g(z)=\int\limits_{[0,z+\Delta z]}T_{B}(z+\Delta z-s)(A-B)T_{A}(s)ds-\int\limits_{[0,z]}T_{B}(z-s)(A-B)T_{A}(s)ds
$$
$$
=\int\limits_{[0,z]}T_{B}(z+\Delta z-s)(A-B)T_{A}(s)ds-\int\limits_{[0,z]}T_{B}(z-s)(A-B)T_{A}(s)ds
$$
$$
+\int\limits_{[z,z+\Delta z]}T_{B}(z+\Delta z-s)(A-B)T_{A}(s)ds
$$
$$
=\left(T_{B}\left(\frac{z}{2}+\Delta z\right)-T_{B}\left(\frac{z}{2}\right)\right)\int\limits_{[0,z]}T_{B}\left(\frac{z}{2}-s\right)(A-B)T_{A}(s)ds
$$
$$
+\int\limits_{[z,z+\Delta z]}T_{B}(z+\Delta z-s)(A-B)T_{A}(s)ds.
$$
It follows that
$$
\|g(z+\Delta z)-g(z)\|_{\mathfrak{S}_1}\leq \left\|T_{B}\left(\frac{z}{2}+\Delta z\right)-T_{B}\left(\frac{z}{2}\right)\right\|M^2\|A-B\|_{\mathfrak{S}_1}
$$
$$
+M^2\|A-B\|_{\mathfrak{S}_1}|\Delta z|\to 0\quad  (\Delta z\to 0).
$$
Therefore the function $f(z):=\mathrm{tr}F(z)$ is bounded and continuous, too. Moreover, since $F$ is analytic in the right half plane with respect to the operator norm, we have for every closed path $C$ located at this half plane that
$$
\oint\limits_C f(z)dz=\mathrm{tr}\oint\limits_C F(z)dz=0.
$$
So by the Morera's Theorem $f$ is  analytic in the right half plane, as well.
By the well known theorem of L. Schwartz there is a unique tempered distribution $\xi$ supported in $\mathbb{R}_+$  such that $f(z)=\langle\xi(t),e^{-zt}\rangle$, the Laplace transform of $\xi$.
Since $\|T_A(u)-T_B(u)\|_{\mathfrak{S}_1}\leq \mathrm{const}\cdot u$ and $\int_{(0,\infty)}ud\mu(u)=\psi^\prime(-0)\ne\infty$, the formula (1) implies that
$$
\mathrm{tr}(\psi(A)-\psi(B))=\int\limits_{(0,\infty)}\mathrm{tr}(T_A(u)-T_B(u))d\mu(u)
$$
$$
=\int\limits_{(0,\infty)}f(u)ud\mu(u)=\int\limits_{(0,\infty)}\langle\xi(t),e^{-ut}\rangle ud\mu(u),
$$
and the proof is complete.

Let us consider the following space of test functions
$$
\mathcal{S}_{\mathbb{R}_+}:=\{\varphi\in C^\infty(\mathbb{R}_+): \|\varphi\|_{j,m}:=\sup\limits_{t\in \mathbb{R}_+}|t^m\varphi^{(j)}(t)|<\infty \forall j, m\in \mathbb{Z}_+\},
$$
endowed with the family of seminorms $\|\cdot\|_{j,m}, j, m\in \mathbb{Z}_+$. It is known that the dual space $(\mathcal{S}_{\mathbb{R}_+})^\prime$ equals to the space  $\mathcal{S^\prime}(\mathbb{R}_+)$ of tempered distributions supported in $\mathbb{R}_+$ \cite[Chapter I, 1.3]{VDZ}.

\begin{corollary}
If, in addition  to the conditions mentioned in Theorem 8, $\psi^\prime(-t)$ belongs to the  space $\mathcal{S}_{\mathbb{R}_+}$, then
$$
\mathrm{tr}(\psi(A)-\psi(B))=\langle\xi(t), \psi^\prime(-t)\rangle.
$$
\end{corollary}

Proof. For every $u> 0$ consider the function $\exp_{-u}(t):=e^{-ut}$ from $\mathcal{S}_{\mathbb{R}_+}$. We claim that the function $u\mapsto \|\exp_{-u}\|_{j,m}$ is integrable with respect to the measure $ud\mu (u)$  $(j, m\in \mathbb{Z}_+)$. Indeed, since $\|\exp_{-u}\|_{j,m}=(m/e)^mu^{j-m}$ (supremum is reached at the point $t=m/u$), it suffices to prove that all functions $u^l, l\in \mathbb{Z}$, are integrable with respect to the measure $\mu$.
To this end note that since the function
$$
\varphi(t):=\psi^\prime(-t)=\int\limits_{(0,\infty)}e^{-ut}ud\mu(u)
$$
belongs to  $\mathcal{S}_{\mathbb{R}_+}$, we have for all  $j, m\in \mathbb{Z}_+$
$$
\|\varphi\|_{j,m}=\sup\limits_{t\in \mathbb{R}_+}\int\limits_{(0,\infty)}t^mu^{j+1}e^{-ut}d\mu(u)<\infty.
$$
If we choose the sequence of reals $t_k$ such that $t_k\uparrow m/u$, then $t_k^me^{-ut_k}\uparrow (m/e)^mu^{-m} (k\to\infty)$ and for all $ j, m\in \mathbb{Z}_+$
$$
\left(\frac{m}{e}\right)^m\int\limits_{(0,\infty)}u^{j+1-m}d\mu(u)=
\lim\limits_{k\to\infty}\int\limits_{(0,\infty)}t_k^mu^{j+1}e^{-ut_k}d\mu(u)\leq \|\varphi\|_{j,m}.
$$
Therefore
$$
\langle\xi, \varphi\rangle=\langle\xi, \int\limits_{(0,\infty)}\exp_{-u}ud\mu(u)\rangle
$$
$$
=\int\limits_{(0,\infty)}\langle\xi,\exp_{-u}\rangle ud\mu(u)=\mathrm{tr}(\psi(A)-\psi(B))
$$
 by Theorem 8.

\begin{corollary}
If, in addition  to the conditions mentioned in Theorem 8, $\xi$ is a measure, then
$$
\mathrm{tr}(\psi(A)-\psi(B))=\int\limits_{\mathbb{R}_+}\psi^\prime(-t)d\xi(t).
$$
\end{corollary}

It follows from Theorem 8 and Tonelli's Theorem.

     \end{document}